\documentclass{amsart}
\usepackage{epsfig}
\usepackage{amssymb,latexsym,ifthen}
\newtheorem{Def}{Definition}

\newtheorem{Thm}{Theorem}

\newtheorem{Rem}{Remark}
\newtheorem{Prop}{Proposition}

\newtheorem{lemma}{Lemma}
\newtheorem{Problem}{Problem}
\newcommand{\inv}{\operatorname{inv}}

\newcommand{\s}{\scriptscriptstyle}
\newcommand{\p}{\scriptstyle}
\newtheorem{Cor}{Corollary}
\newtheorem{Conj}{Conjecture}
\newcommand{\Grass}{\operatorname{W}}

\newcommand{\Arr}{\operatorname{Arr}}

\begin{document}

\title{Quasi-Commuting Families of Quantum Minors}

\author{Josh Scott}
\maketitle

\begin{center} Department of Mathematics \end{center}
\begin{center} Northeastern University \end{center}
\begin{center} Boston, MA 02115 \end{center}
\begin{center} e-mail: josh @ mystic.math.neu.edu \end{center} 
\begin{center} \today \end{center}

\bigskip
\bigskip
\noindent
In [7] a combinatorial criterion for quasi-commutativity is established for pairs of quantum Pl\"ucker coordinates in the
quantized coordinate algebra $\Bbb{C}_q[\mathcal{F}]$ of the flag variety of type $A$. This paper
attempts to generalize these results by producing necessary and sufficient conditions for pairs
of quantums minors in the quantized coordinate algebra $\Bbb{C}_q[Mat_{k \times m}]$ to quasi-commute.
In addition we study the combinatorics of maximal (by inclusion) families of pairwise
quasi-commuting quantum minors and pose relevant conjectures.  \\

\section{Introduction}

\noindent
Let $\Bbb{C}_q[Mat_{k \times m}]$ be the $q$-deformation of the coordinate ring of the space of $k \times m$
complex matrices where $k \leq m$. This is the $\Bbb{C}(q)$-algebra with unity generated
by indeterminates $x_{i,j}$ for $i \in [1 \dots k]$ and
$j \in [1 \dots m]$
subject to the Faddeev-Reshetikhin-Takhtadzhyan relations [2]:
\begin{equation}
\begin{alignat*}{2}
& x_{s,t}x_{i,j} = q \ x_{i,j}x_{s,t} \qquad \qquad &&\mbox{if either} \ s>i \ \
\mbox{and} \ t=j \\
&                                                 && \mbox{\bf or} \ \ \ \ \ \ \ \  s=i \ \
\mbox{and} \ t>j \\
& x_{s,t}x_{i,j} = x_{i,j}x_{s,t}   \qquad \qquad &&\mbox{if} \ s>i \ \ \mbox{and} \ \ t<j \\
& x_{s,t}x_{i,j} = x_{i,j}x_{s,t} + (q - q^{-1}) \ x_{i,t}x_{s,j} \qquad \qquad &&\mbox{if}
\ s>i \ \ \mbox{and} \ \ t>j \end{alignat*}
\end{equation}

\bigskip
\noindent
In this paper we shall be concerned with
a special family of elements $\Delta_{I,J} \in \Bbb{C}_q[Mat_{k \times m}] $ indexed by pairs of non-empty
subsets $I$ and $J$ of
$[1 \dots k]$ and $[1 \dots m]$ respectively with $|I| = |J|= l$. They are defined by:

\[ \Delta_{I,J} := \sum_{\sigma \in S_l} (-q)^{-l(\sigma)} x_{i_1,j_{\sigma(1)}}\cdots
   x_{i_l,j_{\sigma(l)}} \ ,\]

\noindent
where $I = \{ i_1 < \cdots < i_l \}$, $J = \{ j_1 < \cdots < j_l \}$, and $ l(\sigma)$ is the
length of the $l$-permutation $\sigma$. The element $\Delta_{I,J}$ is the $q$-deformation of the classical
determinant and for this reason we call the $\Delta_{I,J}$'s quantum minors.

\begin{Def} Two quantums minors $\Delta_{A,B}$ and $\Delta_{C,D}$ quasi-commute
if $\Delta_{C,D}\Delta_{A,B} \newline = q^c \ \Delta_{A,B}\Delta_{C,D} $ for some integer
$c$. The integer $c$ is uniquely determined by $\Delta_{A,B}$ and $\Delta_{C,D}$ and we will denote
its value by the symbol $ c \bigl( \Delta_{A,B} \  | \ \Delta_{C,D} \bigr) $. Note that
$ c \bigl( \Delta_{C,D} \ | \ \Delta_{A,B} \bigr) = -c \bigl( \Delta_{A,B} \ | \ \Delta_{C,D}
\bigr)$ for any quasi-commuting pair.
\end{Def}

\noindent
We can now state the central problems we will address is this
paper, namely:

\begin{Problem} Find necessary and sufficient conditions
for two quantum minors $\Delta_{A,B}$ and $\Delta_{C,D}$ to quasi-commute.
In addition, explicitly compute $c \bigl( \Delta_{A,B} \ | \ \Delta_{C,D} \bigr) $ in terms
of $A$, $B$, $C$, and $D$.

\end{Problem}

\begin{Problem} Find a combinatorial mechanism which will describe and produce all maximal (by inclusion)
families of pairwise quasi-commuting quantum minors.
\end{Problem}

\noindent
Problems 1 and 2 are motivated by the study of dual canonical bases for quantum groups of type $A$.
It is conjectured in [1], and partially proved in [8],
that products of quasi-commuting quantum minors constitute a part of the
dual canonical basis for the quantum group $\Bbb{C}_q[GL(n,\Bbb{C})]$.
Problem 2 is also motivated by the study of total positivity as described in [3] and [4]. \\

\noindent
Problem 1 is resolved using techniques developed in [7].
Ostensibly Problem 1 is more general than its counterpart in [7]
which only addresses the quantum flag variety. Nevertheless
we demonstrate in this paper that Problem 1 can be reduced to a special case
of the problem treated in [7] - namely the problem of determining when two
quantum Pl\"ucker coordinates of the corresponding quantum Grassmannian quasi-commute.
The criterion for quasi-commutativity is described in
terms of the notion of "weak separability" as put forth in [7].

\begin{Def} Given two subsets $I$ and $J$ of $[1 \dots n]$ we write $I \prec J$ if $i<j$ for all $i \in I$
and all $j \in J$. We say $I$ and $J$ are weakly separated if at least one of the
following two conditions holds: \\

\qquad 1. $|I| \geq |J|$ and $J-I$ can be partitioned into a disjoint union $J-I = J' \sqcup J''$

\indent \indent \indent \quad so that $J' \prec I-J \prec J''$. \\

\qquad 2. $|J| \geq |I|$ and $I-J$ can be partitioned into a disjoint union $I-J = I' \sqcup I''$

\indent \indent \indent \quad so that $I' \prec J-I \prec I''$. \\

\end{Def}

\noindent
We associate to any pair of subsets $A \subset [1 \dots k]$ and
$B \subset [1 \dots m]$ of equal size the subset
$ S(A,B) \subset [1 \dots k+m]$ of size $k$ defined as follows:
\[ S(A,B) =  \Bigl\{ \ b + k \ \Bigl\| \ b \in B \ \Bigr\}  \sqcup [1 \dots k] - w_{\s{0}}(A) \  \]
where  $w_{\s{0}}$ is the order reversing permutation of $[1 \dots k]$.
Problem 1 is settled by the following two Theorems: \\

\begin{Thm}
The quantum minors $\Delta_{A,B}$ and $\Delta_{C,D}$ in $\Bbb{C}_q[Mat_{k,m}]$
quasi-commute if and only if $S(A,B)$ and
$S(C,D)$ are weakly separated subsets of $[1 \dots m+k]$. \\

\end{Thm}

\begin{Thm} Suppose $I = S(A,B)$ and $J = S(C,D)$ are
weakly separated subsets of $[1 \dots m+k]$ satisfying case 1 in Definition 2. Then

\[ c \bigl( \Delta_{A,B} \ | \ \Delta_{C,D} \bigr)  = |J''| - |J'| + |A| - |C|. \] \\

\end{Thm}

\noindent
In proving Theorems 1 and 2 we use a quantum analogue of
the well known embedding of $Mat_{k \times m}$ as an affine chart in the Grassmannian
$\Bbb{G}_{k,k+m}$; this embedding sends a $k \times m$ matrix $(x_{i,j})$ to the row
space of the $k \times (k+m)$ matrix

\setcounter{MaxMatrixCols}{20} \[
\begin{pmatrix} 0 & & & & &  1 & x_{\s{1,1}} & \cdot & \cdot & \cdot & \cdot & x_{\s{1,m}} \\
                  & & & & -1 & & \cdot & & & & & \cdot \\
                  & & & \cdot & & & \cdot & & & & & \cdot \\
                  & & \cdot & & & & \cdot & & & & & \cdot \\
                  & \cdot & & & & & \cdot & & & & & \cdot \\
          \ (-1)^{k-1} & & & & & 0 & x_{\s{k,1}} & \cdot & \cdot & \cdot & \cdot & x_{\s{k,m}} \end{pmatrix}
\]
\setcounter{MaxMatrixCols}{10}

\noindent
The corresponding quantum analogue is an embedding of $\Bbb{C}_q[Mat_{k \times m}]$ into the
quantized coordinate ring $\Bbb{C}_q[\Bbb{G}_{k,k+m}]$ - the so called
quantum Grassmannian as defined in [10]. This embedding allows us to reduce questions about
quantum minors to corresponding questions about quantum Pl\"ucker coordinates. \\

\noindent
Theorem 1 implies that $ \mathcal{C} = \{ \Delta_{A_1,B_1}, \dots, \Delta_{A_s,B_s} \}$ is a maximal
collection of pairwise quasi-commuting quantum minors in $\Bbb{C}_q[Mat_{k,m}]$ if and only if
$\big\{ S(A_1,B_1), \dots , \newline S(A_s,B_s) \big\} \sqcup \big\{ [1 \dots k ] \big\}$ is a maximal collection of
pairwise weakly separated $k$-subsets of $[ 1 \dots k+m]$. This identification is a central component
in our attempt to resolve Problem 2. Theorem 1.3 of [7] asserts that the size of any
maximal collection of pairwise weakly separated $k$-subsets of $[1 \dots n]$ is sharply bounded
by $k(n-k) + 1$. Setting $n = k + m$ we obtain: \\

\begin{Prop}
The size of any maximal collection of pairwise
quasi-commuting quantum minors in $\Bbb{C}_q[Mat_{k \times m}]$ is sharply bounded by $km$. \\
\end{Prop}

\noindent
In [7] the following {\it purity} property is conjectured: all maximal
collections of pairwise weakly separated subsets (not neccessarily $k$-subsets) of $[1 \dots n]$
have size $\binom {n+1} 2 +1$. The analogue of this purity conjecture for $k$-subsets is given by: \\

\begin{Conj}[Purity]
All maximal collections of pairwise weakly separated $k$-subsets
of $[1 \dots n]$ have size $k(n-k) +1$. Equivalently, all maximal
collections of pairwise quasi-commuting quantum minors in
$\Bbb{C}_q[Mat_{k \times m}]$ have size $km$. \\
\end{Conj}

\noindent
In Sections 5 and 6 we prove this assertion for the cases
$k=2$ and $k=3$ respectively. \\

\noindent
In Section 3 we expose a new feature specific to the quantum
Grassmannian: quasi-commutativity of
the quantum Pl\"ucker coordinates in $\Bbb{C}_q[\Bbb{G}_{k,n}]$ is preserved under the natural action of the dihedral
group $D_n$. More precisely, we show that the natural $D_n$-action on $k$-subsets of $[1 \dots n]$
preserves weak separbility.
We do not know of an analogue of this action for the full quantum flag variety.
Let $\Grass(k,n)$ be the set of all maximal collections of pairwise weakly separated $k$-subsets in
$[1 \dots n]$. The induced $D_n$-action on $\Grass(k,n)$ is instrumental in proving several
assertions in this paper. \\

\noindent
For a set $I$ and elements $x$ and $y$ let $I \p{xy} $ denote $I \cup \{x,y\}$.
The set $\Grass(k,n)$ possesses the following interesting structure. \\

\begin{Thm} \label{movethm}
Let $\mathcal{C}$ be a maximal collection of pairwise weakly separated $k$-subsets of $[1 \dots n]$.
Suppose that $I \p{ij}$, $I \p{it}$, $I \p{js}$, $I {\p{st}} \in \mathcal{C}$ for some
$i<s<j<t$ and for some $I \subset [1 \dots n] - \{i,j,s,t\} $ with $|I| = k-2$.
Then $\mathcal{C}$ contains either $I \p{ij}$ or $I\p{st}$ and not both.
Moreover, the transformation

\begin{equation} \mathcal{C} \longmapsto \begin{cases}
                 \mathcal{C} - \{{I \p{ij}}\} \sqcup \{{I \p{st}}\} &  \text{if} \ \ {I \p{ij}} \in \mathcal{C} \\
                 \mathcal{C} - \{{I \p{st}}\} \sqcup \{{I \p{ij}}\} &  \text{if} \ \ {I \p{st}} \in \mathcal{C} \end{cases}
\end{equation} \\

\noindent
preserves weak separability and maximality. \\

\end{Thm}

\noindent
This transformation is an analogue of the Yang-Baxter ``flip" introduced in [7];
here we refer to these transformations as $(2,4)$-moves due to the fact that they originate on
$\Bbb{C}_q[\Bbb{G}_{2,4}]$. \\

\begin{Conj}[Transitivity]
Let $\mathcal{C}$ and $\mathcal{B}$ be any collections in $\Grass(k,n)$. Then there is a sequence of
$(2,4)$-moves transforming \ $\mathcal{C}$ into $\mathcal{B}$. \\
\end{Conj}

\noindent
If true the conjecture effectively settles Problem 2.
In addition it provides a method to obtain all collections in $\Grass(k,n)$: simply propagate
a given maximal collection by all possible $(2,4)$-moves.
In Section 3 we explain why the validity of Conjecture 2 implies the validity of
Conjecture 1. In Sections 5 and 6 we prove this Conjecture 2 for the cases
$k=2$ and $k=3$. In Section 8 we explore applications of this conjecture to total positivity. \\

\noindent
In Section 4 we describe certain maximal collections in $\Grass(k,n)$ arrising from {\it double wiring arrangements}.
In Section 7 we present a construction that recursively generates all collections in $\Grass(3,n)$ by
lifting collections from $\Grass(3,n-1)$. In principle this construction should provide
a method to compute the size of $\Grass(3,n)$.  \\

\section{ The Quantum Grassmannian and Proofs of Theorems 1 and 2}

\begin{Def}

The quantum Grassmannian $\Bbb{C}_q [ \Bbb{G}_{k,n} ]$, as defined in [10], is the $\Bbb{C}(q)$-algebra with unity
generated by all quantum Pl\"ucker coordinates
$\Delta^K$ where $K$ is a $k$-subset of $[1 \dots n]$ subject to the relations:

\[ \sum_{i \in I-J} \ (-q)^{\inv(i,I) - \inv(i,J)} \ \Delta^{I-\{i\}} \ \Delta^{J \sqcup \{i\}} \ = \ 0 \]
\\

\noindent
for any $(k+1)$-subset $I$ and $(k-1)$-subset $J$. Here $\inv(i,X)$ is the number of $x \in X$
such that $i>x$.

\end{Def}

\begin{Prop}[Quantum Stieffel-Pl\"ucker Correspondence]
 There exists a unique $\Bbb{C}(q)$-algebra embedding $\varphi: \Bbb{C}_q[Mat_{k \times m}]
 \longrightarrow \Bbb{C}_q[\Bbb{G}_{k,k+m}]$ such that

 \[ \Delta_{I,J} \longmapsto  q^{\binom{l}{2}} \ \Delta^{l-1} \
    \Delta^{S(I,J)} \] \\

 \noindent
 where $l = |I| =|J|$ and $\Delta = \Delta^{[1 \dots k]}$.

\end{Prop}
\begin{proof}
The proof that the Faddeev-Reshetikhin-Takhtadzhyan relations are
preserved under the correspondence $x_{i,j} \longmapsto \ \Delta^{
S(\{i\},\{j\}) } $ and that $\Delta_{I,J}$ is sent to
$q^{\binom{l}{2}} \ \Delta^{l-1} \ \Delta^{S(I,J)}$ is a simple
modification of the proof of the quantum analogue of Bazin's
theorem presented in Theorem 3.8 of [6]. \\

\noindent The classical analogue of $\varphi$, obtained by
specializing $q$ to $1$, is easily seen to be injective. This
taken together with Theorem 3.5(c) of [10] and the fact that the
monomials consisting of products of lexicographically ordered
generators $x_{i,j}$ form a basis for $\Bbb{C}_q[Mat_{k \times m}]$
over $\Bbb{C}(q)$ proves injectivity of $\varphi$. \\
\end{proof}

\noindent It is well known that $\Delta^{[1 \dots k]}$ is
quasi-central. Thus Proposition 2 tells us that two quantum minors
$\Delta_{A,B}$ and $\Delta_{C,D}$ will quasi-commute exactly when
the corresponding quantum Pl\"ucker coordinates $\Delta^{S(A,B)}$
and $\Delta^{S(C,D)}$ quasi-commute. In turn, the conditions for
two quantum Pl\"ucker coordinates to quasi-commute are explained
by the following proposition of [7]: \\

\begin{Prop} Two quantum Pl\"ucker coordinates $\Delta^I$ and $\Delta^J$ in \ $ \Bbb{C}_q [ \Bbb{G}_{k,n} ]$
quasi-commute if and only if $I$ and $J$ are weakly separated. If $I$ and $J$ satisfy case 1 of
Definition 2 then $c \bigl( \Delta^I \ | \ \Delta^J \bigr) = |J''| - |J'|$. \\

\end{Prop}

\noindent
Theorem 1 now follows from Propositions 2 and 3. Theorem 2 also follows from Propositions 2 and 3 along with
the fact that $c \bigl( \Delta^{|A|-1} \ | \ \Delta^{S(C,D)} \bigr) =
|C|(|A| -1)$ and $ c \bigl( \ \Delta^{ S(A,B) } \ | \
\Delta^{|C| -1} \bigr) = |A|(1 - |C|) $. \\

\section{Proof of Theorem 3 }

\noindent
It is convenient to visualize a $k$-subset of $[1 \dots n]$ as a
subpolygon of the regular polygon with $n$ vertices labeled counter-clockwise by the indices
$[1 \dots n]$. Represent the dihedral group $D_n$ as the group of symmetries of the $n$-gon.
Clearly $D_n$ acts on the set of $k$-subsets of $[1 \dots n]$ under this realization. \\

\begin{Prop}
If two $k$-subsets $I$ and $J$ of $[1 \dots n]$ are weakly
separated then $g(I)$ and $g(J)$ are weakly separated for any $g \in D_n$. \\
\end{Prop}

\begin{proof}

\noindent
In [7] it is shown that $I$ and $J$ are weakly separated precisely when, after interchanging
$I$ and $J$ if neccessary, either: \\

\qquad a) $|I|<|J|$ and there do not exist three indices $a<b<c$ such that

\indent \indent \indent \quad $I \cap \{a,b,c\} = \{b\}$
and $J \cap \{a,b,c\} = \{a,c\}$ \ {\bf or} \\

\qquad b) $|I|=|J|$ and there do not exist four indices $a<b<c<d$ such that

\indent \indent \indent \quad $I \cap \{a,b,c,d\} =
          \{a,c\}$ and $J \cap \{a,b,c,d\} = \{b,d\}$  \\

\noindent
Part b) above indicates that two $k$-subsets $I$ and $J$
are weakly separated precisely, when viewed as subpolygons, no diagonal of the subpolygon $I$
disjoint from $J$ {\it crosses} a diagonal of $J$ disjoint from $I$. This property is
clearly preserved under any dihedral symmetry of the $n$-gon. \\

\end{proof}

\noindent
A $k$-subset $I$ is called {\it boundary} if it consists of $k$
consecutive indices of the $n$-gon; i.e. any $k$-subset of the form $g([1 \dots k])$
for $g \in D_n$. Since $[1 \dots k]$ is weakly separated with every $k$-subset
it follows that the set of all $k$-boundary subsets is common to every
maximal collection of pairwise weakly separated $k$-subsets. \\

\noindent
{\bf Proof of Theorem 3:} \\

\noindent
To prove the first part of the theorem notice that
since $I \p{ij}$ and $I \p{st}$ are not weakly separated it is clear that both can not
be in $\mathcal{C}$. So we need only demonstrate that one of them
is present in $\mathcal{C}$. Given a $k$-subset $J$ of $ [1 \dots n]$ such
that $J$ is weakly separated from $I \p{is}$, $I \p{sj}$, $I \p{jt}$, $I \p{it}$ and
different from $I \p{ij}$ and $I \p{st}$
we need to show that $J$ is weakly separated from both $I \p{ij}$ and $I \p{st}$. \\

\noindent
Proposition 4 shows that we may reduce the proof to the
case of $t=n$ after suitably translating the collection $\mathcal{C}$ by the {\it dihedral action}.
Assume that $t=n$. Let $J^- = J - \{n\}$.
Since $|J| = k$ and $J$ is different from $I \p{ij}$ and $I \p{st}$,
it follows that $J^- $ is different from both $I \p{ij}$ and $I \p{s}$.
By Lemma 3.2 of [7], $J^-$ is weakly separated from $I \p{is}$, $I \p{sj}$, $I \p{j}$, $I \p{i}$.
By Lemma 5.2 of [7], it follows that $J^-$ is weakly separated from both $I \p{ij}$ and $I \p{s}$ and,
after an easy application of part b) above, that $J$ is weakly separated from both $I \p{ij}$ and
$I \p{sn}$, as claimed. \\

\noindent
The above argument also shows that the transformation (1) preserves weak separability and maximality,
thus concluding the proof of Theorem 3. $\Box$ \\ \\

\noindent
Returning to Conjecture 2, notice that if it is true
and if we can find a collection $\mathcal{A}$ in $\Grass(k,n)$ for which
$|\mathcal{A}| = k(n-k)+1$ then Conjecture 1 will follow. One can easily verify that
the collection $\mathcal{A} = \mathcal{A}_n$ whose non-boundary sets are

\[  \Bigl\{ \ [1 \dots i] \sqcup [j \dots k+j-i-1] \ \Big\| \ 1 \leq i <k \ \text{and}
                         \ i+1<j<n+i-k \ \Bigr\} \] \\
has the desired properties. \\

\section{Wiring Arrangements}

\noindent
In [7] a recursive procedure is described through which all maximal families of pairwise
weakly separated subsets (not neccessarily $k$-subsets) of $[1 \dots n]$ are obtained.
In principle this recursion can be restricted to produce all families in $\Grass(k,n)$.
Nevertheless, the process is not very practical.
In this section we explore a non-recursive combinatorial device which parametrizes a large portion
of the collections in $\Grass(k,n)$. This device is a modification of a construction
in [3]. \\

\noindent
Recall that the symmetric group $S_n$ is generated by the simple reflections $s_i = (i,i+1)$
satisfying the Coxeter relations. A reduced word for an element $g \in S_n$ is
sequence of indices $i_1 ,\dots i_l$ such that $g = s_{i_1} \cdots s_{i_l}$ with $l$ minimal.
For the group $S_k \times S_m$
we will use the indices $[\overline{1}, \dots , \overline{k-1} ]$ to
label the simple reflections corresponding to the $S_k$ component and the
indices $[1 \dots m-1]$ to label the simple reflections for the
$S_m$ component. Under this convention a reduced word for an element $(u,v)
\in S_k \times S_m$ can be identified with a shuffle of a reduced word
for $u$, written with indices in $[\overline{1},\dots,\overline{k-1}]$, and
a reduced word for $v$ written with indices $[1 \dots m-1]$. \\

\noindent
Let $w_{\s{0}}^{\s{(k)}}$ and $w_{\s{0}}^{\s{(m)}}$ denote the longest elements in
$S_k$ and $S_m$ respectively.
We say a reduced word for $(w_{\s{0}}^{\s{(k)}},w_{\s{0}}^{\s{(m)}}) \in S_k \times S_m$ is
{\it optimal} if the associated reduced word for $w_{\s{0}}^{\s{(m)}}\in S_m$ has
a total of only $\binom{m-k}2$ occurrences of the indices $[k+1,\dots m-1]$.
Given an opitimal reduced word ${\bf i}$ of $ (w_{\s{0}}^{\s{(k)}},w_{\s{0}}^{\s{(m)}})
\in S_k \times S_m$ we will manufacture a maximal collection $\mathcal{C}({\bf i})$ of
pairwise quasi-commuting quantum minors. This collection is obtained by means of the {\it double wiring
arrangement} $\Arr({\bf i})$ attached to ${\bf i}$, as introduced in [3]. \\

\noindent
Recall first the definition of a {\it single wiring arrangement} attached to a reduced word.
It is easiest to understand this definition with an example. Consider the
reduced word $1231$ of the permutation $\begin{pmatrix} 1 & 2 & 3 & 4 \\ 3 & 2 & 4 & 1 \end{pmatrix}
\in S_4 $. The corresponding single wiring arrangement is:

\begin{figure}[ht]

\setlength{\unitlength}{1.0pt}

\begin{center}
\begin{picture}(270,60)(6,0)
\thicklines
%\dark{

  %\setlength{\linethickness}{2pt}

  \put(0,0){\line(1,0){40}}
  \put(40,0){\line(1,1){20}}
  \put(60,20){\line(1,0){40}}
  \put(100,20){\line(1,1){20}}
  \put(120,40){\line(1,0){40}}
  \put(160,40){\line(1,1){20}}
  \put(180,60){\line(1,0){100}}

  \put(0,20){\line(1,0){40}}
  \put(40,20){\line(1,-1){20}}
        \put(60,0){\line(1,0){160}}
  \put(220,0){\line(1,1){20}}
  \put(240,20){\line(1,0){40}}

  \put(0,40){\line(1,0){100}}
  \put(100,40){\line(1,-1){20}}
  \put(120,20){\line(1,0){100}}
  \put(220,20){\line(1,-1){20}}
  \put(240,0){\line(1,0){40}}

  \put(0,60){\line(1,0){160}}
  \put(160,60){\line(1,-1){20}}
  \put(180,40){\line(1,0){100}}

  \put(50,-10){\bf 1}
  \put(110,-10){\bf 2}
  \put(170,-10){\bf 3}
  \put(230,-10){\bf 1}

%}

\end{picture}
\end{center}
\caption{Single wiring arrangement}
\label{fig:double-wiring}
\end{figure}
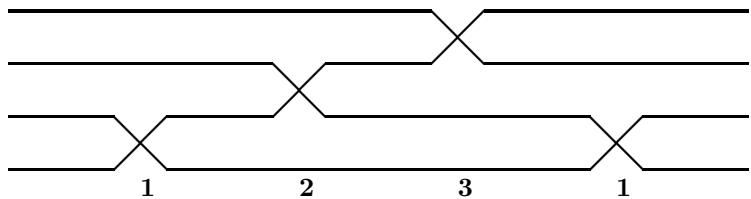

\noindent
We associate a crossing at the $i$th level (counting from the bottom up) for each $i$ in the reduced word.
To obtain the double wiring arrangement for $(u,v) \in S_k \times S_m$ we superimpose the
single wiring arrangements for the reduced words of $u$ and $v$ respectively aligning them
closely in the vertical direction (starting at the bottom) and intertwining their respective crossings
as dictated by the shuffle. To distinguish the two wiring arrangements we colour the diagram for $u$ red.
For example, the double wiring arrangement corresponding to the reduced word

\[  {\bf i} = 2 \ \overline{1} \ 1 \
2 \ 3 \ \overline{2} \ 2 \ 1 \ 4 \ \overline{1} \ 3 \ 2 \ 1 \ \ \text{for} \ \
(w_{\s{0}}^{\s{(3)}},w_{\s{0}}^{\s{(5)}}) \in S_3 \times S_5 \ \ \text{is:} \]

\begin{figure}[ht]

\setlength{\unitlength}{.5 pt}

\begin{center}

\begin{picture}(640,160)(6,0)

\thicklines

%\dark{

  %\setlength{\linethickness}{2pt}

  \put(0,0){\line(1,0){100}}
  \put(100,0){\line(1,2){20}}
  \put(120,40){\line(1,0){40}}
  \put(160,40){\line(1,2){20}}
  \put(180,80){\line(1,0){40}}
  \put(220,80){\line(1,2){20}}
  \put(240,120){\line(1,0){160}}
  \put(400,120){\line(1,2){20}}
  \put(420,160){\line(1,0){220}}

  \put(0,40){\line(1,0){40}}
  \put(40,40){\line(1,2){20}}
  \put(60,80){\line(1,0){100}}
  \put(160,80){\line(1,-2){20}}
  \put(180,40){\line(1,0){120}}
  \put(300,40){\line(1,2){20}}
  \put(320,80){\line(1,0){140}}
  \put(460,80){\line(1,2){20}}
  \put(480,120){\line(1,0){160}}

  \put(0,80){\line(1,0){40}}
  \put(40,80){\line(1,-2){20}}
  \put(60,40){\line(1,0){40}}
  \put(100,40){\line(1,-2){20}}
  \put(120,0){\line(1,0){220}}
  \put(340,0){\line(1,2){20}}
  \put(360,40){\line(1,0){160}}
  \put(520,40){\line(1,2){20}}
  \put(540,80){\line(1,0){100}}

  \put(0,120){\line(1,0){220}}
  \put(220,120){\line(1,-2){20}}
  \put(240,80){\line(1,0){60}}
  \put(300,80){\line(1,-2){20}}
  \put(320,40){\line(1,0){20}}
  \put(340,40){\line(1,-2){20}}
  \put(360,0){\line(1,0){220}}
  \put(580,0){\line(1,2){20}}
  \put(600,40){\line(1,0){40}}

  \put(0,160){\line(1,0){400}}
  \put(400,160){\line(1,-2){20}}
  \put(420,120){\line(1,0){40}}
  \put(460,120){\line(1,-2){20}}
  \put(480,80){\line(1,0){40}}
  \put(520,80){\line(1,-2){20}}
  \put(540,40){\line(1,0){40}}
  \put(580,40){\line(1,-2){20}}
  \put(600,0){\line(1,0){40}}

  \put(40,-25){\bf 2}
  \put(65,-25){$\bar{1}$}
  \put(105,-25){\bf 1}
  \put(163,-25){\bf 2}
  \put(223,-25){\bf 3}
  \put(245,-25){$\bar{2}$}
  \put(303,-25){\bf 2}
  \put(345,-25){\bf 1}
  \put(405,-25){\bf 4}
  \put(425,-25){$\bar{1}$}
  \put(465,-25){\bf 3}
  \put(525,-25){\bf 2}
  \put(585,-25){\bf 1}

%}

%\light{ \thinlines

  \put(0,5){\line(1,0){60}}
  \put(60,5){\line(1,2){20}}
  \put(80,45){\line(1,0){160}}
  \put(240,45){\line(1,2){20}}
  \put(260,85){\line(1,0){380}}

  \put(0,45){\line(1,0){60}}
  \put(60,45){\line(1,-2){20}}
  \put(80,5){\line(1,0){340}}
  \put(420,5){\line(1,2){20}}
  \put(440,45){\line(1,0){200}}

  \put(0,85){\line(1,0){240}}
  \put(240,85){\line(1,-2){20}}
  \put(260,45){\line(1,0){160}}
  \put(420,45){\line(1,-2){20}}
  \put(440,5){\line(1,0){200}}

%}

\end{picture}

\end{center}

\caption{Double wiring arrangement}

\label{fig:double-wiring}

\end{figure}
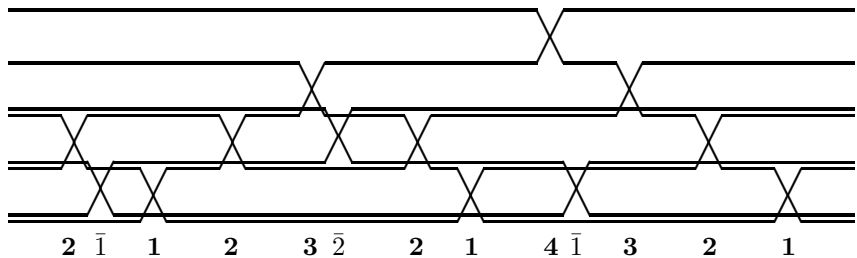

\noindent
To obtain the collection $\mathcal{C}({\bf i})$ label the black wires $1$ through $m$ bottom-up at the
left hand side of the arrangement and label the red wires $1$ through $k$ bottom-up at
the right hand side of the arrangement. Label each chamber $C$
in the first $k$ strips of the arrangement with $I(C)$ - the set of labels
of red lines passing beneath the chamber - and $J(C)$ - the set of black line labels passing beneath the chamber.
For example the above double wiring arrangement is labeled

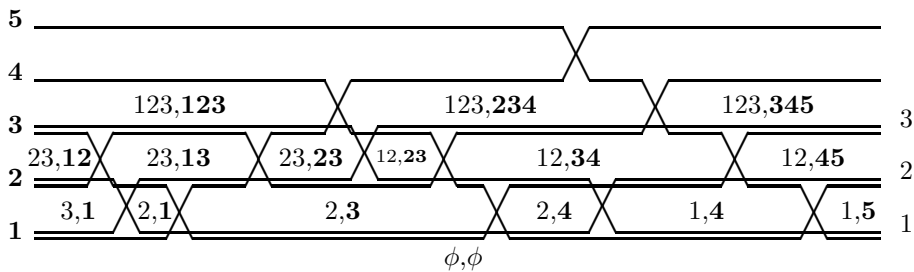
\begin{figure}[ht]

\setlength{\unitlength}{.5 pt}

\begin{center}

\begin{picture}(640,170)(6,0)

\thicklines

%\dark{

  %\setlength{\linethickness}{2pt}

  \put(0,0){\line(1,0){100}}
  \put(100,0){\line(1,2){20}}
  \put(120,40){\line(1,0){40}}
  \put(160,40){\line(1,2){20}}
  \put(180,80){\line(1,0){40}}
  \put(220,80){\line(1,2){20}}
  \put(240,120){\line(1,0){160}}
  \put(400,120){\line(1,2){20}}
  \put(420,160){\line(1,0){220}}

  \put(0,40){\line(1,0){40}}
  \put(40,40){\line(1,2){20}}
  \put(60,80){\line(1,0){100}}
  \put(160,80){\line(1,-2){20}}
  \put(180,40){\line(1,0){120}}
  \put(300,40){\line(1,2){20}}
  \put(320,80){\line(1,0){140}}
  \put(460,80){\line(1,2){20}}
  \put(480,120){\line(1,0){160}}

  \put(0,80){\line(1,0){40}}
  \put(40,80){\line(1,-2){20}}
  \put(60,40){\line(1,0){40}}
  \put(100,40){\line(1,-2){20}}
  \put(120,0){\line(1,0){220}}
  \put(340,0){\line(1,2){20}}
  \put(360,40){\line(1,0){160}}
  \put(520,40){\line(1,2){20}}
  \put(540,80){\line(1,0){100}}

  \put(0,120){\line(1,0){220}}
  \put(220,120){\line(1,-2){20}}
  \put(240,80){\line(1,0){60}}
  \put(300,80){\line(1,-2){20}}
  \put(320,40){\line(1,0){20}}
  \put(340,40){\line(1,-2){20}}
  \put(360,0){\line(1,0){220}}
  \put(580,0){\line(1,2){20}}
  \put(600,40){\line(1,0){40}}

  \put(0,160){\line(1,0){400}}
  \put(400,160){\line(1,-2){20}}
  \put(420,120){\line(1,0){40}}
  \put(460,120){\line(1,-2){20}}
  \put(480,80){\line(1,0){40}}
  \put(520,80){\line(1,-2){20}}
  \put(540,40){\line(1,0){40}}
  \put(580,40){\line(1,-2){20}}
  \put(600,0){\line(1,0){40}}

  \put(-20,0){\bf 1}
  \put(-20,40){\bf 2}
  \put(-20,80){\bf 3}
  \put(-20,120){\bf 4}
  \put(-20,160){\bf 5}
%}

%\light{ \thinlines

  \put(0,5){\line(1,0){60}}
  \put(60,5){\line(1,2){20}}
  \put(80,45){\line(1,0){160}}
  \put(240,45){\line(1,2){20}}
  \put(260,85){\line(1,0){380}}

  \put(0,45){\line(1,0){60}}
  \put(60,45){\line(1,-2){20}}
  \put(80,5){\line(1,0){340}}
  \put(420,5){\line(1,2){20}}
  \put(440,45){\line(1,0){200}}

  \put(0,85){\line(1,0){240}}
  \put(240,85){\line(1,-2){20}}
  \put(260,45){\line(1,0){160}}
  \put(420,45){\line(1,-2){20}}
  \put(440,5){\line(1,0){200}}

  \put(655,5){1}
  \put(655,45){2}
  \put(655,85){3}
%}

  \put(20,15){3,{\bf 1}}
  \put(78,15){2,{\bf 1}}
  \put(220,15){2,{\bf 3}}
  \put(380,15){2,{\bf 4}}
  \put(495,15){1,{\bf 4}}
  \put(610,15){1,{\bf 5}}

  \put(-5,55){23,{\bf 12}}
  \put(85,55){23,{\bf 13}}
  \put(185,55){23,{\bf 23}}
  \put(252,58){ $\scriptstyle{12,{\bf 23}}$ }
  \put(380,55){12,{\bf 34}}
  \put(565,55){12,{\bf 45}}

  \put(75,95){123,{\bf 123}}
  \put(310,95){123,{\bf 234}}
  \put(520,95){123,{\bf 345}}

  \put(310,-20){$\phi$,${\bf \phi}$}

\end{picture}

\end{center}

\caption{Labeled arrangement}

\label{fig:double-wiring}

\end{figure}

\bigskip
\bigskip
\bigskip
\noindent
Let $\mathcal{C}({\bf i}) = \Bigl\{ \ \Delta_{I(C),J(C)} \ \Bigl\| \ C \
\text{a chamber of $\Arr({\bf i})$ of level $\leq k$ } \ \Bigr\}$. \\

\begin{lemma}
Let ${\bf i}$ be an optimal reduced word for $(w_{\s{0}}^{\s{(k)}},w_{\s{0}}^{\s{(m)}}) \in S_k \times S_m$. Then
the size of $\mathcal{C}({\bf i})$ is $km$.
\end{lemma}

\begin{proof}
Given ${\bf i}$, the number of chambers in the first $k$ strips of the corresponding double wiring
arrangement is equal to the number of red and black crossings in the first $k$ strips plus $k$ - corresponding
to the $k$ far right chambers . The number of black (respectively red) crossings in
the first $k$ strips in turn is given by the number of simple reflections $j$ (respectively $\bar{j}$)
occurring in the reduced word ${\bf i}$ with $1 \leq j \leq k$. The number of $\bar{j}$ in ${\bf i}$
with $1 \leq j \leq k$ is $\binom k 2$. The number of of $j$ in ${\bf i}$ with $1 \leq j \leq k$
is $\binom m 2  - \# \ \big\{ \ j \ \text{occurring in ${\bf i}$} \ \big\| \ k+1 \leq j \leq m-1 \ \big\}$;
if ${\bf i}$ is {\it optimal} this will be $\binom m 2  - \binom {m-k} 2$.
Consequently the number of chambers occurring in the first $k$ strips of the double wiring arrangement
for ${\bf i}$ optimal ( or equivalently the size of $\mathcal{C}({\bf i})$ ) is:

\[ \binom k 2 + \binom m 2   - \binom {m-k} 2 + k  = mk   \] \\

\end{proof}

\begin{Prop}
If ${\bf i}$ is an optimal reduced word for $(w_{\s{0}}^{\s{(k)}},w_{\s{0}}^{\s{(m)}}) \in S_k \times S_m$ then
$\mathcal{C}({\bf i})$ is a maximal collection of pairwise quasi-commuting
quantum minors in $\Bbb{C}_q[Mat_{k \times m}]$. Moreover,
given $\Delta_{A,B}$ and $\Delta_{I,J}$ in $\mathcal{C}({\bf i})$ either

\begin{equation} A - I \prec I - A \ \ \text{and} \ \ J - B \prec B - J \ \ \text{\bf or} \end{equation}
\begin{equation} I - A \prec A - I \ \ \text{and} \ \ B - J \prec J - B \ \ \text{\ \   } \end{equation} \\

\end{Prop}

\begin{proof}
Take any quantum minors $\Delta_{A,B}$ and $\Delta_{I,J}$ in $\mathcal{C}({\bf i})$. Lemma 4.1
of [7] proves that if $X$ and $Y$ are chamber sets of a single wiring arrangement then
either $X-Y \prec Y-X$ or $Y-X \prec X-Y$. This, taken together with the fact that
the single wiring arrangements for the $S_k$ and $S_m$ components of ${\bf i}$ are oppositely
labeled, proves the second part of the proposition. \\

\noindent
To prove that $\Delta_{A,B}$ and $\Delta_{I,J}$ quasi-commute we must show that
$S(A,B)$ and $S(I,J)$ are weakly separated. We may assume, after exchanging $A$ with $I$ and
$B$ with $J$ if neccessary, that $A-I \prec I-A$ and $J-B \prec B-J$. This in turn is
equivalent to

\[ \Big( S(A,B) - S(I,J) \Big) \cap [1 \dots k] \prec
   S(I,J) - S(A,B) \prec \Big( S(A,B) - S(I,J) \Big) - [1 \dots k]  \] \\

\noindent
which demonstrates that $S(A,B)$ and $S(I,J)$ are weakly separated. The fact that
$\mathcal{C}({\bf i})$ is maximal follows from Lemma 1 and Proposition 1. \\
\end{proof}

\noindent
It is possible to prove the converse of Proposition 5, namely: If $\mathcal{C}$ is
a collection of quantum minors $\Delta_{A,B}$ whose indices
pairwise satisfy either condition 2 or 3, and if $\mathcal{C}$ is maximal with respect to this
property, then $\mathcal{C}$ is of the form $\mathcal{C}({\bf i})$ for some
optimal reduced word ${\bf i}$. \\

\noindent
Given an optimal reduced word ${\bf i}$ the following collection is in $\Grass(k,k+m)$:

\[ \Big\{ \ S\Big(I(C),J(C)\Big) \ \Big\| \ C \ \text{a chamber of} \Arr({\bf i}) \ \text{of level}
\ \leq k \ \Big\} \sqcup \Big\{ [1 \dots k] \Big\} \] \\

\noindent
In the case of $\Grass(3,6)$ all collections are obtained via double wiring arrangements.
There are $34$ in total and they are explicitly described
in [3] and [4]. Every maximal family in $\Grass(3,6)$ is dihedrally equivalent
to one of the following five collections (we omit boundary sets):

\[ \Bigl\{ \ \{124\}, \{125\}, \{134\}, \{145\} \ \Bigr\} \ \ \ \
\Bigl\{ \ \{124\}, \{125\}, \{145\}, \{245\} \ \Bigr\} \]

\[ \Bigl\{ \ \{124\}, \{134\}, \{145\}, \{146\} \ \Bigr\} \ \ \ \
\Bigl\{ \ \{125\}, \{134\}, \{135\}, \{145\} \ \Bigr\} \]

\[ \Bigl\{ \ \{135\}, \{136\}, \{145\}, \{235\} \ \Bigr\} \] \\

\noindent
In general it is not the case that every maximal collection in $\Grass(k,n)$
corresponds to some double wiring arrangement, even after dihedral
translation. This is evidenced already in the case of $\Bbb{C}_q[\Bbb{G}_{2,n}]$.
In Section 5 we shall demonstrate such a maximal collection. \\

\section{The case of $\Bbb{C}_q[\Bbb{G}_{2,n}]$}

\noindent
%Label the vertices of an $n$-gon with indices $[1 \dots n] $ in a counter-clockwise fashion.
We identify the 2-subsets of $[1 \dots n]$ with chords inscribed in a regular $n$-gon.
Clearly two 2-subsets of $[1 \dots n]$ are weakly separated
if and only if the corresponding chords do not cross in the interior of the polygon.
Under this identification collections $\mathcal{C} \in \Grass(2,n)$ correspond to maximal collections of non-crossing
chords - i.e. triangulations of an $n$-gon.

\begin{Thm}[Transitivity] \label{transthm1}

Let $\mathcal{C},\mathcal{B} \in \Grass(2,n)$. Then there is a sequence of $(2,4)$-moves
transforming  $\mathcal{C}$ into $\mathcal{B}$.
\end{Thm}

\begin{proof}
This theorem follows from the well known fact that the any two triangulations are connected
by a series of chord exchanges where the diagonal chord of an inscribed quadralateral is "flipped" to its crossing pair.
The diagonal "flips" correspond to $(2,4)$-moves.
\end{proof}

\begin{Cor}[Purity] \label{puritythm1}
Let $\mathcal{C} \in \Grass(2,n)$. Then $|C| = 2(n-2)+1$.
\end{Cor}

\begin{proof}
Immediate corollary of Theorem \ref{transthm1}. \\
\end{proof}

\noindent
Since $\Grass(2,n)$ is identified with the set of triangulations of an $n$-gon it follows that
$| \Grass(2,n) | $ is the Catalan number
${1 \over {n-1}}\binom{2n-4}{n-2}$. For $k>2$ the size of $\Grass(k,n)$ is not known. \\

\noindent
In [5] it is shown that the coordinate ring $\Bbb{C}[\Bbb{G}_{2 \times n}]$ has a basis consisting of
all monomials of Pl\"ucker coordinates whose indices are
%lexicographically ordered and
pairwise weakly separated.
Using the quantum short Pl\"ucker relation given by

\[ \Delta^{I \p{ij}} \ \Delta^{I \p{st}} \ = \ q \ \Delta^{I \p{is}} \ \Delta^{I \p{jt}} \ + \ q^{-1} \ \Delta^{I \p{it}} \ \Delta^{I \p{sj}} \]

\noindent
for $i<s<j<t$ as a straightening rule, we obtain the following quantum analogue of this result:

\begin{Prop}
The set of all monomials consisting of lexicographically ordered pairwise quasi-commuting quantum Pl\"ucker coordinates is a basis
for $\Bbb{C}_q[\Bbb{G}_{2,n}]$. \\
\end{Prop}

\noindent
Using Proposition 5 and the identification of maximal collections in $\Grass(2,n)$ with triangulations
of an $n$-gon we can characterize those maximal collections which can be parametrized, up to
the dihedral action, by double wiring arrangements.
Given $\mathcal{C} \in \Grass(2,n)$ there exists $g \in D_n$ for which $g \cdot \mathcal{C}$ is
parametrized by a double wiring arrangement if and only if there exists an external edge of
the polygon (i.e. a boundary $2$-set) such that for any other external edge there is {\bf no} chord
in the associated triangulation, which separates both the edges and is disjoint from both.
The following collection in $\Grass(2,9)$, represented as a triangulation,  is an example of a collection
which is not parametrized, up to the dihedral action, by a double wiring arrangement: \\

\begin{figure}[ht]

\setlength{\unitlength}{.5 pt}

\begin{center}
\begin{picture}(320,280)(6,0)
\includegraphics{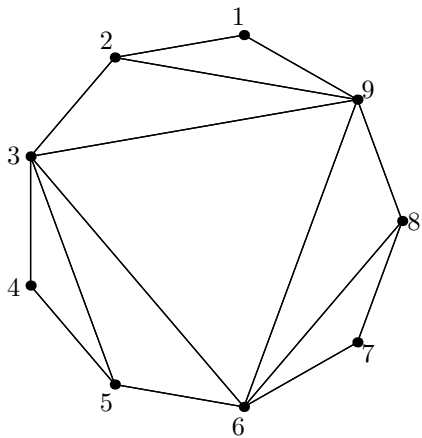}

\put(-145,305){1}
\put(-245,287){2}
\put(-315,200){3}
\put(-315,100){4}
\put(-245,13){5}
\put(-145,-5){6}
\put(-47,50){7}
\put(-47,250){9}
\put(-12,150){8}

\end{picture}
\end{center}

\caption{Non-Parametrized $\Grass(2,9)$ collection} 
\end{figure}

\section{The case of $\Bbb{C}_q[\Bbb{G}_{3,n}]$}

\noindent
In this section we prove the Transitivity and Purity Conjectures
for $k=3$. \\

\begin{Thm}[Transitivity] \label{transthm2}
Let $\mathcal{C},\mathcal{B} \in \Grass(3,n)$. Then there is a sequence of $(2,4)$-moves
transforming  $\mathcal{C}$ into $\mathcal{B}$. \\
\end{Thm}

\begin{Cor}[Purity] \label{puritythm2}
Let $\mathcal{C} \in \Grass(3,n)$ then $|\mathcal{C}| = 3(n-3)+1$. \\
\end{Cor}

\noindent
{\bf Proof of Transitivity:} \\

\noindent
The essential strategy is to show that any collection $\mathcal{C} \in \Grass(3,n)$ can be reduced
by a sequence of $(2,4)$-moves to the "base" collection
$\mathcal{A}_n$ whose non-boundary $3$-sets are

\[ \Bigl\{ \ \{1,s,s+1\} \ \Big\| \ 2<s<n-1 \ \Bigr\} \ \bigsqcup \ \
   \Bigl\{ \ \{1,2,s\} \ \Big\| \ 3<s<n \ \Bigr\} \] \\

\noindent
We first prove that whenever a collection
$\mathcal{C}$ can be $(2,4)$-reduced to $\mathcal{A}_n$ then so can any of its dihedral translations
$g \cdot \mathcal{C}$ for $g \in D_n$. In Lemma \ref{dihedralprop}
we then show that any maximal collection can be translated dihedrally to a maximal collection
containing the $3$-set $\{1,n-2,n-1\}$.
%This restricts the focus of the conjecture
%to collections which contain the $3$-set $\{1,n-2,n-1\}$.
We conclude the proof by showing that any such collection can
be reduced by a sequence of $(2,4)$-moves to the collection $\mathcal{A}_n$.   \\

\begin{lemma} \label{translemma1}
Let $\mathcal{C} \in \Grass(3,n)$. If $\mathcal{C}$ can be reduced by a sequence of $(2,4)$-moves
to $\mathcal{A}_n$ then so can the collection $g \cdot \mathcal{C}$ for any $g \in D_n$.
\end{lemma}

\begin{proof}
Since the $D_n$-action preserves $(2,4)$-moves
it is enough to verify this assertion in the case where $\mathcal{C} = \mathcal{A}_n$. \\

\noindent
Proceed by induction on $n$. For $n \leq 4$ the statement is evident. Assume $n>4$.
It is enough to verify the claim for the group elements $\rho_n$ and $\sigma_n$, which
generate $D_n$, given by

\[ \rho_n = \begin{pmatrix} 1 & 2 & \cdots & n-1 & n \\ 2 & 3 & \cdots & n & 1 \end{pmatrix} \quad
\sigma_n = \begin{pmatrix} 1 & 2 & 3 & 4 & 5 & \cdots \\ 2 & 1 & n & n-1 & n-2 & \cdots \end{pmatrix} \] \\

\noindent
This follows from the observation that if $g \cdot
\mathcal{C}$ can be reduced by a sequence of $(2,4)$-moves to $\mathcal{B}$ then $hg \cdot \mathcal{C}$ can
be reduced to $h \cdot \mathcal{B}$. \\

\noindent
The collection $\sigma_n \cdot \mathcal{A}_n$ contains the $3$-sets $\{1,2,n-1\}$, $\{2,n-2,n-1\}$,
$\{n-2,n-1,n\}$, $\{1,n-1,n\}$, and $\{2,n-1,n\}$. Applying the $(2,4)$-move which
replaces $\{2,n-1,n\}$ with $\{1,n-2,n-1\}$ we obtain $\sigma_{\s{n-1}} \cdot \mathcal{A}_{n-1}
\sqcup \Bigl\{ \ \{1,2,n\} , \{1,n-1,n\}, \{n-2,n-1,n\} \ \Bigr\} $. By induction
$\sigma_{\s{n-1}} \cdot \mathcal{A}_{n-1}$ can be reduced by a sequence of $(2,4)$-moves to $\mathcal{A}_{n-1}$.
Thus $\sigma_n \cdot \mathcal{A}_n$ can be reduced to $\mathcal{A}_{n-1} \sqcup \Bigl\{ \ \{1,2,n\}, \{1,n-1,n\},
\{n-2,n-1,n\} \ \Bigr\} = \mathcal{A}_n$. \\

\noindent
To deal with $\rho_n$, notice that $\rho_n \cdot \mathcal{A}_n$ contains the 3-sets $\{1,2,n\}$,
$\{ 1,2,3\}$, $\{2,3,n-1\}$, $\{2,n-1,n\}$, and $\{2,3,n\}$. We apply the $(2,4)$-move
which replaces $\{2,3,n\}$ with $\{1,2,n-1\}$. This new collection contains the 3-sets $\{1,n-1,n\}$ ,
$\{1,2,n-1\}$, $\{2,n-2,n-1\}$, $\{n-2,n-1,n\}$, and $\{2,n-1,n\}$. We may apply the $(2,4)$-move
which replaces $\{2,n-1,n\}$ with $\{1,n-1,n-2\}$. The resulting collection is exactly $\rho_{\s{n-1}}
\cdot \mathcal{A}_{n-1} \sqcup \Bigl\{ \ \{1,2,n\} , \{1,n-1,n\}, \{n-2,n-1,n\} \ \Bigr\}$. By
the induction hypothesis $\rho_{\s{n-1}} \cdot \mathcal{A}_{n-1}$ can be reduced by a sequence
of $(2,4)$-moves to $\mathcal{A}_{n-1}$. Consequenlty $\rho_n \cdot \mathcal{A}_n$ can
be reduced to $\mathcal{A}_{n-1} \sqcup \Bigl\{ \ \{1,2,n\}, \{1,n-1,n\}, \{n-2,n-1,n\} \ \Bigr\} =
\mathcal{A}_n$. \\

\end{proof}

\begin{lemma} \label{dihedralprop}
Given $\mathcal{C} \in \Grass(3,n)$ there exists $g \in D_n$ such that $g \cdot \mathcal{C}$ contains the
$3$-set $\{1,n-2,n-1\}$.
\end{lemma}

\begin{proof}
For a $3$-subset $I$ of $[1 \dots n]$ define the diameter of $I$
to be the minimal cardinality of a boundary $k$-subset of $[1 \dots n]$ that contains $I$.
Thus the boundary $3$-subsets are precisely those of diameter $3$.
Let us call $3$-subsets of diameter $4$ {\it almost boundary} subsets.
It suffices to prove that every maximal collection $\mathcal{C}$ contains an almost boundary subset.\\

\noindent
Assume by contradiction that $\mathcal{C}$ does not contain an almost boundary $3$-subset.
We make the following easy observation:

\begin{Rem}
Let $a$, $b$, $c$, and $d$ be four consecutive vertices in $[1 \dots n]$;
then the $3$-subsets that are not weakly separated with an almost boundary subset
$\{a,c,d\}$ are precisely the non-boundary $3$-subsets containing $b$ but
not $a$.
\end{Rem}

\noindent
Therefore our assumption and maximality of $\mathcal{C}$ imply
that for every two consecutive vertices $a$ and
$b$ in $[1 \dots n]$, there is a non-boundary $3$-subset in $\mathcal{C}$ which
contains $b$ but not $a$. \\

\noindent
Choose a non-boundary $3$-subset $\{a, c, d \}$ in $\mathcal{C}$ of minimal possible
diameter.
Without loss of generality, we can assume that a
boundary subset of minimal cardinality that contains $\{a,c,d \}$
has $a$ and $d$ as its endpoints; let us denote
this boundary subset by $[a,d]$.
We can also assume that $c$ is not a neighbor of $a$.
Let $b$ be the neighbor of $a$ in $[a,d]$.
Consider a $3$-subset $I$ in $\mathcal{C}$ such that
$I$ contains $b$ but not $a$.
Since $I$ is weakly separated from $\{a,c,d \}$ it must be contained in
$[b,d] = [a,d] - \{a\}$.
But then $I$ has smaller diameter than $\{a,c,d \}$ which contradicts
our choice of $\{a,c,d \}$.
This proves the claim and hence the lemma as well. \\

\end{proof}

\noindent
For any collection $\mathcal{C} \in \Grass(3,n)$ we define its {\it height} $H(\mathcal{C})$ to be the number
of non-boundary $3$-sets containing $n$.
An immediate consequence of Remark 1 is that $H(\mathcal{C}) = 0$ if and only if both
$\{1,2,n-1\}$ and $\{1,n-2,n-1\}$ are in $\mathcal{C}$.

\begin{lemma} \label{translemma2}
Let $\mathcal{C} \in \Grass(3,n)$ with $\{1,n-2,n-1\} \in \mathcal{C}$. Then $\mathcal{C}$ can be reduced by
a sequence of $(2,4)$-moves to a collection of height $H = 0$.
\end{lemma}

\begin{proof}
We proceed by induction on the height.
If $H(\mathcal{C})=0$ then we are already done. Assume inductively that the assertion is true for collections of height $H = k \geq 0$
and let $\mathcal{C}$ be a collection of height $H(\mathcal{C}) = k+1$.
We need the following:  \\

\begin{lemma} \label{pinchlemma} Let $\mathcal{C} \in \Grass(3,n)$ and suppose that $\{1,n-2,n-1\} \in \mathcal{C}$. Then
there exists a unique index $b>1$ such that both $\{1,b,n-1\}$ and $\{1,b,n\}$ are in $\mathcal{C}$.
We call $b$ the {\it pinch point} over $n$ and $n-1$.
\end{lemma}

\begin{proof}
Let $b$ be the maximal index with the property that $\{1,b,n\} \in \mathcal{C}$. Suppose, by
contradiction, that $\{1,b,n-1\} \notin \mathcal{C}$.
By maximality of $\mathcal{C}$ this means there exists a non-boundary set $I \in \mathcal{C}$ which is not weakly
separated with $\{1,b,n-1\}$. Therefore there exist indices $s,t \in I$ such that
one of the following holds: \\

\qquad \qquad \qquad \qquad \qquad 1. \ $1<s<b<t<n-1$

\qquad \qquad \qquad \qquad \qquad 2. \ $1<s<b$ and $t=n$

\qquad \qquad \qquad \qquad \qquad 3. \ $b<s<n-1$ and $t=n$ \\

\noindent
Case 1: Since $I$ and $\{1,b,n\}$ are weakly separated it follows that $b \in I$. But then $I$ will be weakly separated with $\{1,b,n-1\}$. \\

\noindent
Case 2: Since $\{1,n-2,n-1\} \in \mathcal{C}$ and since $I$ is a non-boundary set containing $n$ it follows that $1 \in I$. But
then $I$ will be weakly separated with $\{1,b,n-1\}$. \\

\noindent
Case 3: Once again it must be the case that $1 \in I$. So $I = \{1,s,n\}$ where $b<s$ violating the maximalitly of $b$. \\

\noindent
Hence $\{1,b,n-1\} \in \mathcal{C}$. Suppose there was another pinch point $b' \ne b$.
Either $b'<b$ or $b'>b$. If $b'<b$ then $\{1,b',n-1\}$ will not be weakly separated from $\{1,b,n\}$.
If $b'>b$ then $\{1,b',n\}$ will not be weakly separated from $\{1,b,n-1\}$. Both possibilities
violate that fact that $\mathcal{C}$ consists of only pairwise weakly separated $3$-sets.
Uniqueness follows.

\end{proof}

\begin{lemma} \label{keyprop2}
Let $\mathcal{C} \in \Grass(3,n)$ and assume $\{1,n-2,n-1\} \in \mathcal{C}$. Let $b$ be the pinch
point over $n$ and $n-1$. Assume in addition that $b>2$.
Then there exists $a$ with $1<a<b$ such that both
$\{1,a,b\}$ and $\{1,a,n\}$ are in $\mathcal{C}$.
\end{lemma}

\begin{proof}
Consider the set of all $x$ with the property that $x<b$ and $\{1,x,n\} \in \mathcal{C}$.
This set is clearly non-empty since $2<b$ and $\{1,2,n\} \in \mathcal{C}$.
Let $a$ be the maximal index with this property. Suppose $\{1,a,b\} \notin \mathcal{C}$.
Then there exists $I \in \mathcal{C}$ with $s,t \in I$ such that one of the following holds: \\

\qquad \qquad \qquad \qquad \qquad 1. \ $1<s<a<t<b$

\qquad \qquad \qquad \qquad \qquad 2. \ $1<s<a<b<t$

\qquad \qquad \qquad \qquad \qquad 3. \ $a<s<b<t$ \\

\noindent
Case 1: Since $\{1,a,n\} \in \mathcal{C}$ it follows that $I$ and $\{1,a,n\}$ must be weakly separated.
The only way this can happen is that $a \in I$. But then $I$ and $\{1,a,b\}$ will be weakly separated. \\

\noindent
Case 2: Since $I$ and $\{1,a,n\}$ are weakly separated it must be the case that $t = n$.
Since $\{1,n-2,n-1\} \in \mathcal{C}$ it follows that $I$ and $\{1,n-2,n-1\}$ are weakly separated.
The only way this can be resolved is that $1 \in I$. But then $I$ and $\{1,a,b\}$ are weakly separated. \\

\noindent
Case 3: Either $t=n$ or not. Suppose $t \ne n$. Since $\{1,b,n\} \in \mathcal{C}$, and
hence weakly separated from $I$, it follows that $b \in I$ in which
case $I$ and $\{1,a,b\}$ will be weakly separated. Thus $t=n$. Since $\{1,b,n-1\} \in \mathcal{C}$ we
know that $I$ and $\{1,b,n-1\}$ are weakly separated. The only way this can happen is that $1 \in I$ and
hence $I = \{1,s,n\}$. But this violates the maximality of $a$ since $a<s<b$. \\

\noindent
Thus $\{1,a,b\}$ and $\{1,a,n\}$ are in $\mathcal{C}$ as required. \\

\end{proof}

\noindent
Returning to Lemma \ref{translemma2}, let $b$ be the pinch point
of $\mathcal{C}$ - i.e. the unique index $b$ such that both $\{1,b,n-1\}$ and $\{1,b,n\}$ are in $\mathcal{C}$.
If $b=2$ it follows that $\{1,2,n-1\} \in \mathcal{C}$.
This, taken together with the fact that $\{1,n-2,n-1\} \in \mathcal{C}$, violates the hypothesis that $H(\mathcal{C}) > 0$.
Therefore $b>2$. \\

\noindent
Since $b>2$ Lemma \ref{keyprop2} implies that there exists $a$ with $1<a<b$ such that both $\{1,a,b\}$ and
$\{1,a,n\}$ are in $\mathcal{C}$. Thus $\mathcal{C}$ contains $\{1,a,b\}$, $\{1,a,n\}$, $\{1,b,n-1\}$, $\{1,b,n\}$,
and $\{1,n-1,n\}$. The associated $(2,4)$-move for this quintuple replaces $\{1,b,n\}$ with $\{1,a,n-1\}$.  Let
$\mathcal{B}$ be the resulting collection. Notice that $\mathcal{B}$ contains $\{1,n-2,n-1\}$ and that $H(\mathcal{B}) =
H(\mathcal{C}) - 1 = k$. By induction $\mathcal{B}$ can be further reduced by a sequence of $(2,4)$-moves into
a collection of height $H=0$. Concantenating this $(2,4)$-reduction with the $(2,4)$-move transforming $\mathcal{C}$
to $\mathcal{B}$ we obtain the desired reduction for $\mathcal{C}$. \\

\end{proof}

\noindent Now we are ready to finish the proof of Transitivity.
%proceed by induction on $n$.
%The case of $n \leq 4$ is trivial. Assume the assertion if holds for $n \geq 4$.
Let $\mathcal{C} \in \Grass(3,n)$. By Lemma \ref{dihedralprop} there is $g \in D_n$ such that
$g \cdot \mathcal{C}$ contains the $3$-set $\{1,n-2,n-1\}$. By Lemma \ref{translemma2} the collection $g \cdot
\mathcal{C}$ can be reduced by a sequence of $(2,4)$-moves to a collection $\mathcal{B}$ with height $H(\mathcal{B}) = 0$.
The collection $\mathcal{B} - \Bigl\{ \ \{1,2,n\}, \{1,n-1,n\}, \{n-2,n-1,n\} \ \Bigr\}$ is in $\Grass(3,n-1)$ and by
induction on $n$ we can assume that it can be reduced by a sequence of $(2,4)$-moves to $\mathcal{A}_{n-1}$.
Equivalently $\mathcal{B}$ can
be reduced by a sequence of $(2,4)$-moves to $\mathcal{A}_n$. Consequently $g \cdot \mathcal{C}$ can be reduced
to $\mathcal{A}_n$ and applying Lemma \ref{translemma1} we conclude that $\mathcal{C}$ can be reduced to $\mathcal{A}_n$
as required. \\

\section{Reduction}
\noindent
In this section we present a recursive procedure to generate collections in $\Grass(3,n)$. \\

\noindent
Given a $3$-subset $I$ of $[1 \dots n]$, we define

\[ I'  = \begin{cases}
              I \sqcup \{n-1\} - \{n\} &  \text{if} \ n \in I \ \text{and} \ n-1 \notin I \\
                                   \phi&  \text{if} \ n \in I \ \text{and} \ n-1 \in I \\
               I&  \text{if} \ n \notin I  \end{cases} \] \\

\noindent
For $\mathcal{C} \in \Grass(3,n)$ let $\mathcal{C}' = \{ \ I' \ | \ I \in \mathcal{C} \ \}$,
and define $F_{\s{\mathcal{C}}}$ to be
the set of indices $b \in [2 \dots n-1]$ with $\{1,b,n\} \in \mathcal{C}$ such that
$\{1,b\} - \{s,t\} \prec \{s,t\} - \{1,b\} $ whenever $\{s,t,n\} \in \mathcal{C}$ for $1<s<t$.
If  $\mathcal{C}$ contains $\{1,n-2,n-1\}$, let $b_{\s{\mathcal{C}}}$ be the pinch point of $\mathcal{C}$
(see Lemma~\ref{pinchlemma}), that is, the unique index such that both
$ \{1,b_{\s{\mathcal{C}}},n-1\}$ and $\{1,b_{\s{\mathcal{C}}},n\}$ are in $\mathcal{C}$.

\begin{Thm}[Reduction] \label{Redthm} Let $n \geq 4$.
The mapping \ $\mathcal{C} \longmapsto \big( \mathcal{C}' , b_{\s{\mathcal{C}}} \big)$ \ defines a
bijection between collections in $\Grass(3,n)$ containing $\{1,n-2,n-1\}$ and the set

\[ %\mathcal{S} \ :=
\ \Bigr\{ \ \big( \mathcal{B} , b \big) \in \Grass(3,n-1) \times [2 \dots n-2] \
\Big\| \ b \in F_{\s{\mathcal{B}}} \ \Bigr\} \] \\

\noindent
The inverse bijection sends a pair $\big( \mathcal{B}, b \big) $ to the collection $ \hat{\mathcal{B}}_b :=
\big\{ \ I_b \ \big\| \ I \in \mathcal{B} \ \big\} \sqcup \Big\{  \{1,b,n-1\}, %\newline
\{1,n-1,n\}, \{n-2,n-1,n\}  \Big\}$ where

\[ I_b = \begin{cases}
              I - \{n-1\} \sqcup \{n\}&  \text{if} \ n-1 \in I \ \text{and} \ I- \{1,b,n-1\} \prec \ \{1,b\} - I \\
              I&  \text{otherwise}  \end{cases} \] \\

\end{Thm}

\noindent
Since by Lemma~\ref{dihedralprop}, every collection in $\Grass(3,n)$ is dihedrally equivalent to
one containing the near boundary subset
$\{1,n-2,n-1\}$, it follows from Theorem \ref{Redthm} that all collections in $\Grass(3,n)$
can be obtained by first lifting collections in $\Grass(3,n-1)$ by the inverse of the
reduction procedure and then translating them suitably by the dihedral action. \\

\noindent
{\bf Proof of Reduction Theorem:} \\

\noindent
The following lemma shows that
the mapping $\mathcal{C} \longmapsto \big( \mathcal{C}' , b_{\s{\mathcal{C}}} \big)$
is well defined.

\begin{lemma} \label{blue1}
Let $\mathcal{C} \in \Grass(3,n)$. Then $\mathcal{C}' \in \Grass(3,n-1)$, and
$b_{\s{\mathcal{C}}} \in F_{\s{\mathcal{C}'}}$.
\end{lemma}

\begin{proof}
Momentary consideration reveals that $\mathcal{C}'$ consists of pairwise weakly separated $3$-subsets of $[1 \dots n-1]$.
In virtue of Corollary \ref{puritythm2} we know that $\mathcal{C}'$ will be maximal if and only if $|\mathcal{C}'| =
3(n-4)+1$. Since $\{1,n-2,n-1\} \in \mathcal{C}$ it follows that if $I \in \mathcal{C}$ and $I' = \phi$ then either
$I = \{1,n-1,n\}$ or $I = \{n-2,n-1,n\}$. Consequently $|\mathcal{C}'| \leq |\mathcal{C}| -2$.
For $I,J \in \mathcal{C}$ if $I' = J'$ then either $I=J$ or else there exists $b \in [2 \dots n-2]$ such that,
after interchanging $I$ and $J$ if neccessary, $I = \{1,b,n-1\}$ and $J= \{1,b,n\}$. By Lemma \ref{pinchlemma},
$b$ is unique. Hence $|\mathcal{C}'| = |\mathcal{C}| - 3 = 3(n-4) + 1$ as required.
The inclusion $b_{\s{\mathcal{C}}} \in F_{\s{\mathcal{C}'}}$ is
also clear from the definitions.

\end{proof}

\noindent
To prove that the inverse correspondence is well defined, we need to show that $\hat{\mathcal{B}}_b
\in \Grass(3,n)$ and $\{1,n-2,n-1\} \in \hat{\mathcal{B}}_b$
for any $\mathcal{B} \in \Grass(3,n-1)$ and  $b \in F_{\s{\mathcal{B}}} $.
Simple consideration shows that all $3$-subsets in $\hat{\mathcal{B}}_b$ are weakly separated
because $b \in F_{\s{\mathcal{B}}}$.
Since $\mathcal{B}$ is maximal we know by Corollary \ref{puritythm2} that $|\mathcal{B}| = 3(n-4) +1$ and thus
$|\hat{\mathcal{B}}_b| = |\mathcal{B}| + 3 = 3(n-3) + 1$. Corollary \ref{puritythm2} implies that $\hat{\mathcal{B}}_b
\in \Grass(3,n)$. Notice also that $\{1,n-2,n-1\} \in \hat{\mathcal{B}}_b$ since $b \leq n-2$. \\

\noindent
It remains to show that the mappings
$\mathcal{C} \longmapsto \big( \mathcal{C}' , b_{\s{\mathcal{C}}} \big)$
and $\big( \mathcal{B}, b \big) \longmapsto  \hat{\mathcal{B}}_b$
are inverse to each other.
First suppose that $\mathcal{C} = \hat{\mathcal{B}}_b$.
Since both $\{1,b,n\}$ and $\{1,b,n-1\}$ are in $\hat{\mathcal{B}}_b$,
the desired equality $\big( \mathcal{C}' , b_{\s{\mathcal{C}}} \big) =  \big( \mathcal{B} , b \big)$
follows from Lemma \ref{pinchlemma}.
Finally, the equality $\widehat{\mathcal{C}_b'} = \mathcal{C}$ for $b=b_{\s{\mathcal{C}}}$
is clear from the definitions. $\Box$ \\

\noindent
{\bf Example:}
Let $\mathcal{C}$ be the collection in $\Grass(3,6)$ whose non-boundary $3$-sets are

\[ \Big\{ \{ 136\} , \{146\}, \{236\}, \{346\} \Big\} \] \\

\noindent
Here $F_{\s{\mathcal{C}}} = \{ 2, 3\}$. Notice that $4 \notin F_{\s{\mathcal{C}}}$ because
$\{1,4\} - \{23\} \not \prec \{2,3\} - \{1,4\}$. The index $5$ is not present
for the same reason. The two possible lifts of $\mathcal{C}$ (omitting boundaries) are:

\[ \hat{\mathcal{C}}_2 = \Big\{ \{126\}, \{136\}, \{146\}, \{156\}, \{236\}, \{346\} \Big\} \]
\[ \hat{\mathcal{C}}_3 = \Big\{ \{137\}, \{136\}, \{146\}, \{156\}, \{236\}, \{346\} \Big\} \] \\

\section{Positivity}

\noindent
Let $\Bbb{G}_{k,n}(\Bbb{C})$ be the Grassmannian of $k$-subspaces in $\Bbb{C}^n$.
Recall that any $k$-subspace in $\Bbb{G}_{k,n}(\Bbb{C})$ can be represented by a
$k \times n$ matrix whose rows span the $k$-subspace. The Pl\"ucker coordinates
are the maximal minors of this $k \times n$ matrix. We say
a point $p \in \Bbb{G}_{k,n}(\Bbb{C})$ is positive if
it can be represented by a $k \times n$ matrix whose Pl\"ucker coordinates $\Delta^I(p)$
are positive real numbers.

\begin{Def} Let $\mathcal{C}$ be a collection of $k$-subsets of $[1 \dots n]$. We say
that $\mathcal{C}$ is a positivity test if $p \in \Bbb{G}_{k,n}(\Bbb{C})$ is positive if and only if all
$\Delta^I(p)$ are real and positive for each $I \in \mathcal{C}$.

\end{Def}

\noindent
In [7] it is conjectured that maximal families of
pairwise weakly separated subsets (not necessarily $k$-subsets)
of $[1 \dots n]$ give rise to positivity tests for the flag variety
of type $A_n$. The analogue of this result for the
Grassmannian $\Bbb{G}_{k,n}(\Bbb{C})$ is:

\begin{Thm} Let $k=2$ or $k=3$. If $\mathcal{C}$ is a maximal collection of
pairwise weakly separated $k$-subsets of $[1 \dots n]$ then the associated collection of Pl\"ucker
coordinates $\{ \ \Delta^I \ | \ I \in \mathcal{C} \ \}$ is a positivity test.
\end{Thm}

\begin{proof}
Let $\mathcal{C} \in \Grass(k,n)$ and suppose that all $\Delta^I (p)$ are real and positive
for $I \in \mathcal{C}$. We need to show that all other Pl\"ucker coordinates
$\Delta^J (p)$ are real and positive. Take any $J \notin \mathcal{C}$.
Take any maximal collection $\mathcal{B}$ containing $J$. Since $k$ is either
$2$ or $3$ we know that Conjecture 2 holds and thus $\mathcal{C}$ and $\mathcal{B}$ are connected
by a sequence of $(2,4)$-moves. \\

\noindent {\bf Claim:}
Suppose $\mathcal{A}$ is in $\Grass(k,n)$ and is a positivity test. Let $\mathcal{B}$ be in
$\Grass(k,n)$ and assume that $\mathcal{B}$ is obtained from $\mathcal{A}$ by a single $(2,4)$-move.
Then $\mathcal{B}$ is a positivity test. \\

\noindent
Indeed, since $\mathcal{A}$ and $\mathcal{B}$ differ
by a single $(2,4)$-move there exist $i<s<j<t$ and $I$, where $I$ is empty if $k=2$ and $|I|=1$ if $k=3$, such that
$I \p{is}$, $I \p{sj}$, $I \p{jt}$, and $I \p{it}$ are in both $\mathcal{A}$ and $\mathcal{B}$ and such that,
without loss of generality, $\mathcal{B}$ is obtained from $\mathcal{A}$ by replacing $I \p{ij}$ with $I \p{st}$.
The fact that $\mathcal{B}$ is a positivity test is an immediate consequence of the short Pl\"ucker relation

\[ \Delta^{I \p{ij}} \ \Delta^{I \p{st}} \ = \ \Delta^{I \p{is}} \ \Delta^{I \p{jt}} \ + \ \Delta^{I \p{it}} \ \Delta^{I \p{sj}} \]

\noindent
Let $l$ be the minimal number of $(2,4)$-moves required to join $\mathcal{B}$ and $\mathcal{C}$. To prove the theorem
proceed by induction on $l$ and use the claim. \\

\end{proof}

\noindent
A positivity test $\mathcal{C}$ is {\bf minimal} if it has no proper subset which
is also a positivity test. We conjecture that $\mathcal{C}$ is a minimal positivity test
for $\Bbb{G}_{k,n}(\Bbb{C})$ if and only if $\mathcal{C}$ is in $\Grass(k,n)$.
In addition, A. Zelevinsky and S. Fomin conjecture that collections $\mathcal{C}$ in $\Grass(k,n)$ have the
property that any Pl\"ucker coordinate $\Delta^J$ can be uniquely expressed as
a positive Laurent polynomial in the Pl\"ucker coordinates $\Delta^I$ for $I \in \mathcal{C}$.
The author intends to investigate these issues related to positivity in a forthcoming article. \\

\noindent
{\bf Acknowledgements:}
Many thanks to my advisor A. Zelevinsky for numerous detailed conversations and suggestions.
Thanks also to my wife E. Machkasova for help in writing GAP code for computations.

\end{document}